# Dyck Numbers, I. Successor Function

Gennady Eremin

ergenns@gmail.com

October 3, 2022

**Abstract.** Dyck paths are among the most heavily studied Catalan families. In the paper we are dealing with the minimal numbering of Dyck paths, with the resulting numbers, the terms of the OEIS sequence A036991, which we have called Dyck numbers. We consider the successor function on the Dyck numbers, closed formulas are obtained. In particular, the formula for the successor of the Mersenne number is calculated (since OEIS A000225 is a subset of A036991). We have obtained the corresponding algorithm of the Dyck successor function, DS-function, and a Python program has been compiled. At the end of the paper, the size of the A036991 ranges is associated with the terms of OEIS A001405. The corresponding hypothesis is formulated.

## 1 Introduction

Dyck paths are among the most heavily studied Catalan families [1]. In the paper we are dealing with the numbering of Dyck paths, with the resulting numbers, which we have called the *Dyck numbers*. In computability theory, numbering is the assignment of natural numbers to a set of objects, such as functions, graphs, words in some formal language, etc. We can say that numbering is a numerical ordering of objects that makes it easier to refer to them.

The Dyck alphabet includes two characters: up step $U = (1, 1)$ and down step $D = (1, -1)$ (in the case of bracket expressions, we are dealing with left and right parentheses, respectively). In Dyck word, the characters of the alphabet are precisely balanced; in such words, *Dyck dynamics* are observed [2]. Each Dyck path has the same number of U's and D's, and this is the first rule of Dyck dynamics. Any up step has a matching down step and they must be correctly nested.

The second rule of Dyck dynamics is as follows: in any initial fragment of Dyck path (aka Dyck word prefix), the number of U's is at least the number of D's. Zoltan Kasa proposed a simple formula for the second rule (see [3], p. 111). Let the Dyck word have length $2n$ (there are $n$ U's and $n$ D's); in this case, the $i$-th down step, $d_i$, satisfies the inequality

$$2i \leq d_i \leq n + i, \ 1 \leq i \leq n.$$

We get the Dyck path numbering if we simply encode the words in the binary alphabet {0, 1}. The On-Line Encyclopedia of Integer Sequences (OEIS, [4]) has a sequence of binary numbers A063171, which is obtained by encoding Dyck words in the following way: the up step is replaced by 1, and the down step is replaced by zero. Here is the be-ginning of the A063171 (terms sorted in ascending order, initial term 0 corresponds to the empty Dyck path):

> 0, 10, 1010, 1100, 101010, 101100, 110010, 110100, 111000, 10101010 …

As you can see, all terms are even (multiples of 2, binary suffix 0). The OEIS sequence A014486 shows the decimal representations of the binary codes:

> 0, 2, 10, 12, 42, 44, 50, 52, 56, 170, 172, 178, 180, 184, 202, 204, 210 …

In this paper, we will be interested in the reverse binary encoding of Dyck words, the OEIS sequence A036991. In this sequence, we get much smaller natural numbers (something like the *minimal numbering*), all terms are odd, among which there are many prime numbers.

We consider a significant group (bush) of OEIS sequences associated with A036991. Two sequences of such a bush are developed by the author, these are A350346 and A350577.

## 2 Minimal numbering of Dyck paths

In conventional decimal expansion, the reverse Dyck path (or Dyck word) codes are shown in OEIS A036991 in ascending order. Here is the beginning of the A036991:

> 0, 1, 3, 5, 7, 11, 13, 15, 19, 21, 23, 27, 29, 31, 39, 43, 45, 47, 51, 53, 55 …

The reverse binary encoding of Dyck words is OEIS A350346: the up step is replaced by 0, and the down step is replaced by 1. Here are the first A350346 terms (leading zeros are shown in small print):

> 0, 01, 0011, 0101, 000111, 001011, 001101, 00001111, 010011, 010101, 00010111 …

In natural numbers, it is customary to omit leading zeros (as in any number system). Without leading zeros, binary codes are greatly reduced:

(1)    0, 1, 11, 101, 111, 1011, 1101, 1111, 10011, 10101, 10111, 11011, 11101 …

And again, in sequence (1), the initial null term corresponds to an empty Dyck path (or an empty Dyck word).

More often we will work with binary strings. Each non-zero A350346 term starts and ends with unit, so there are fewer zeros than units. Obviously, there is no need for the first rule of Dyck dynamics, since in any binary string it is possible to restore the leading zeros and by recoding back get the corresponding Dyck path. Therefore, we can simplify the Dyck dynamics. It all comes down to a single rule: *in each final fragment (suffix) of the binary code, the number of 0's does not exceed the number of 1's.* Further, working with numerical sequences, we will adhere to exactly this formulation of the Dyck dynamics.

It is easy to see that Kasa's formula will take the following form for the position of the *i*-th zero in the binary code, $z_i$ (recall that in integers, digits are numbered from right to left, starting from 0):

$$2i - 1 \leq z_i \leq n + i - 1, \quad 1 \leq i \leq n.$$

In each binary term, the number of leading zeros varies from 1 to *n* (assuming the corresponding Dyck word has length 2*n*). Without leading zeros the length of such terms varies from *n* to 2*n*-1.

Thus, using OEIS A036991, we get hopefully the minimal numbering of Dyck paths. At the moment, there is no more suitable numbering option. For example, in [5], for each natural number, based on prime factorization, the correct bracket sequence is calculated (among other things, here huge natural numbers correspond to compact Dyck words). In our case, in A036991, each Dyck word corresponds to a certain natural number; but for most natural numbers there are no corresponding Dyck paths. Now we can formulate an obvious definition.

**Definition 1.** A *Dyck number* is a natural number, in the binary expansion of which in each suffix the number of zeros does not exceed the number of units.

The following statement is obvious.

**Proposition 2.** There is a one-to-one correspondence (bijection) A036991 and the set of Dyck paths.

Dyck numbers are fairly common in the OEIS. The author found more than twenty numerical sequences, the terms of which are either Dyck numbers or are associated with the A036991 in some way. We can say that we are dealing with an extensive bush of interconnected sequences. There are hundreds of thousands of sequences in the OEIS, and the author believes that the interested reader can add my list.

Further in the formulas, we will try to adhere to the notation adopted in the OEIS. For example, a(*n*) is the *n*-th term of the current sequence, * is a multiplication, ^ is a power, sqrt is the square root, etc. Binary expansion can be considered as strings (or words), in this case the concatenation operation is convenient, let's denote it ∘. Grouping of repeated fragments of words is possible, for example,

$$1011111_2 = 10 \circ 1^5, \quad 11101101101101_2 = 11 \circ 101 \circ 101 \circ 101 \circ 101 = 11 \circ (101)^4.$$

For many sequences in the OEIS, the following recurrent formula is true:

(2) $$a(n+1) = 2\hat{\ }s * a(n) + t,$$

where *s* is a non-negative integer, and *t* is an arbitrary integer (possibly negative). For example, for odd natural numbers (OEIS A005408), we get a well-known recurrent formula a(*n*+1) = a(*n*) + 2 if we take *s* = 0 and *t* = 2.

In sequence (1), the terms are arranged in order of increasing length of the binary code. A group of terms of the same length is often called a *range*. For example, in the 4-th range there are three binary terms of length 4: 1011, 1101 and 1111.

In each range, the last binary term is a *repunits* (repeating units) [7], a term from OEIS A002275:  0, 1, 11, 111, 1111 … Here the Dyck dynamics is observed, since binary re-punits does not contain zeros at all. In binary expansion, the Dyck word is compressed as much as possible. In decimal representation, binary repunits are Mersenne numbers of the form

$$M_n = 2^n - 1, \ n \geq 0.$$

The first few Mersenne numbers are 0, 1, 3, 7, 15, 31, 63, 127, 255… (OEIS A000225). For Mersenne numbers, the recurrent formula (2) takes the form $a(n+1) = 2*a(n) + 1$, in this case $s = t = 1$. Let's formulate an obvious statement.

**Proposition 3.** Every Mersenne number is a Dyck number, or A000225 ⊂ A036991 and, respectively, A002275 ⊂ A350346.

Mersenne numbers play an important role in A036991. These numbers form the structure of A036991 by breaking the sequence into ranges. The Dyck number $M_n$ is the last term in the $n$-th range, in binary expansion the next term has a length of $n+1$. Let's show the first ranges in A036991/ A350346.

$R_0 = \{M_0\} = \{0\}$,
$R_1 = \{M_1\} = \{1_2\}$,
$R_2 = \{M_2\} = \{11_2\}$,
$R_3 = \{5, M_3\} = \{101_2, 111_2\}$,
$R_4 = \{11, 13, M_4\} = \{1011_2, 1101_2, 1111_2\}$,
$R_5 = \{19, 21, 23, 27, 29, M_5\} = \{10011_2, 10101_2, …, 11101_2, 11111_2\}$ and so on.

Recall that zero term does not correspond to a real Dyck path. Therefore, we will sometimes omit range $R_0$. In fact, many sources do this, for example, Wikipedia or Wolfram MathWorld [6, Sep 2022] ignore the number $M_0 = 0$.

Of interest are the initial terms in the ranges, i.e., these are the terms that immediately follow the Mersenne numbers. Here is the beginning of the sequence of first numbers in the ranges (starting or offset from 1):

(3)     1, 3, 5, 11, 19, 39, 71, 143, 271, 543, 1055, 2111, 4159, 8319, 16511 …

The sequence (3) is not yet in the OEIS. We see some patterns. For example, for the terms 1, 5, 19, 71 … (odd $n = 1, 3, 5, 7$ …), formula (2) is the same as for the Mersenne numbers: $2*a(n) + 1 = a(n+1)$. And for the terms 3, 11, 39, 143 … (even $n = 2, 4, 6$ …), formula (2) is valid in the case $s = 0$ and $t = 2\wedge(n-1)$: $a(n) + 2\wedge(n-1) = a(n+1)$.

In sequence (3) each term is a successor of the corresponding Mersenne number. The successor function considered below will allow us to obtain a closed formula for the terms of (3).

## 3    Successor function on the Dyck numbers

In mathematics, the successor function sends a natural number to the next one: $S(n) = n + 1$. We also need a similar function on the Dyck numbers.

**Definition 4.** The *Dyck Successor function*, *DS-function*, sends any Dyck number to the next one. Let's denote the DS-function by *DS*. For example, $DS(0) = 1$ and $DS(11) = DS(DS(7)) = 13$.

For an arbitrary Dyck number $d$, you can get the successor by iterating over subsequent odd numbers $d+2, d+4$ and so on, checking Dyck dynamics each time. That's usually how it's done. We will offer the reader a fairly compact algorithm and closed formulas for the DS-function.

### 3.1. DS-function for Mersenne numbers.
Let's calculate DS-function for the *n*-th Mersenne numbers, the Dyck number with a binary code of length *n* and without binary zeros. The resulting value will allow us to determine the jump between adjacent ranges. Let's look at the first Mersenne numbers:

$DS(M_0) = 0 + 1 = 0 + 2^0 = M_1 = 1,$
$DS(M_1) = 1 + 2 = 1 + 2^1 = M_2 = 3,$
$DS(M_2) = 3 + 2 = 3 + 2^1 = 5,$
$DS(M_3) = 7 + 4 = 7 + 2^2 = 11,$
$DS(M_4) = 15 + 4 = 15 + 2^2 = 19,$
$DS(M_5) = 31 + 8 = 31 + 2^3 = 39,$
$DS(M_6) = 63 + 8 = 63 + 2^3 = 71,$
$DS(M_7) = 127 + 16 = 127 + 2^4 = 143$ and so on.

The trend is obvious. Further, we will consider non-zero Mersenne numbers. However, the formula (4) considered below is also true for $M_0$. Now let's formulate the following theorem.

**Theorem 5.** For a non-zero n-th Mersenne number

(4) $$DS(M_n) \;=\; M_n \;+ M_m + 1 \;=\; M_n + 2^m, \; m = \lceil n/2 \rceil.$$

*Proof.* We will work with binary expansions of natural numbers. Let's add 1 to the *n*-th Mersenne number: $M_n + 1 = 2^n$. Then in a binary code of length $n+1$, we get 1 and $n$ trailing 0's. Next, in the zero suffix of length $n$, *n*-suffix, we need to place a minimal binary code in which the Dyck dynamics is performed (recall that in each suffix the number of 0's should not exceed the number of 1's). It is logical to fill the right half of the *n*-suffix, *m*-suffix, $m = \lceil n/2 \rceil$, with units. In the case of an odd $n$, there will be a little more 1's. That is, in fact, we put the *m*-th Mersenne number in the *m*-suffix. As a result, we get the required value:

$$2^n + M_m = M_n + 1 + M_m = M_n + 2^m. \qquad \square$$

Now we can write a closed formula for the terms of the sequence (3):

(5) $$a(n) \;=\; DS(M_{n-1}), \; n \geq 1.$$

Obviously, in binary expansion, the formula (4) can be written like this

(6) $$DS(M_n) \;=\; 1 \cdot 0^{n-m} \cdot 1^m \; \text{(in base 2)}.$$

Obviously, the length of the zero fragment is $n - m = n - \lceil n/2 \rceil = \lfloor n/2 \rfloor$. In the future we will show that the Dyck successor function works transparently in the domain of Mersenne numbers.

**3.2. DS-function on the Dyck numbers (not Mersenne numbers).** The binary expansion of any Mersenne number has no internal zeros, so we obtained equality (4) relatively easily. Let's get a verbal algorithm for the successor function for some given Dyck number whose binary expansion has internal zeros (well, at least one zero). But first we will consider the usual Dyck path in the first quadrant of the *xy*-plane.

**Example 6.** Figure 1 shows a slightly modified Dyck path $p$ from [8]. The path returns once to ground level, namely, the valley $v_2$ has zero height (zero $y$ coordinate).

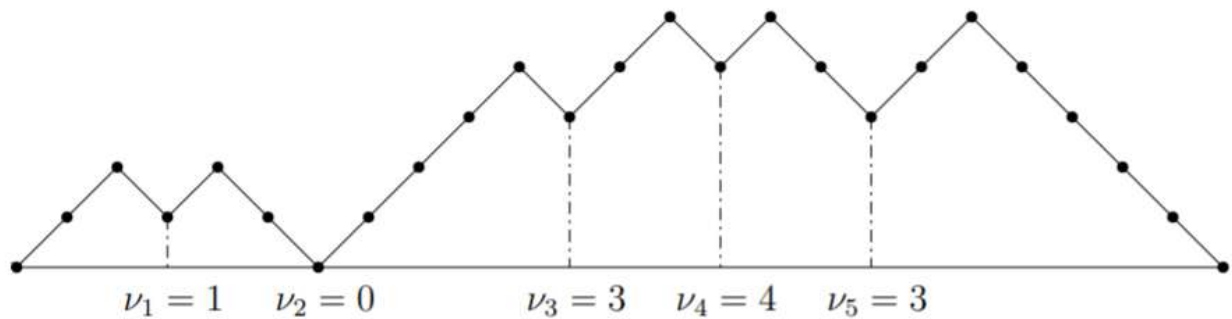

Figure 1. Dyck path of semilength 12.

Consider the various encodings of such a path:

– the string of up steps and down steps

$$p = \mathsf{U^2DUD^2U^4DU^2DUD^2U^2D^5};$$

– the binary expansion of the corresponding Dyck number 2,893,215 is

(7) $\qquad d = \text{00}\,1\,0\,11\,0000\,1\,00\,1\,0\,11\,00\,11111$ (two leading zeros);

– the heights of the vertices, the line is filled from right to left (five valleys are marked in red)

(8) $\qquad h(d) = 2\,1\,21\,0123\,4\,34\,5\,4\,54\,34\,54321.$

For convenience, we will assume that in (8) each digit does not exceed 9 (but we can also work with hexadecimal codes).

Since the Dyck word is encoded as a binary string, and in a binary number the digits are numbered from right to left (and we also remove leading zeros), it is logical to number word attributes from right to left as well. Accordingly, we calculate the heights in $h(d)$ from right to left. It is easy to see that for the binary expansion of the Dyck number, the height of the binary digit is the difference between the number of units and the number of zeros located on the right (in the adjacent suffix). The height of the rightmost unit is 1.

**Definition 7.** Let's call the height of the deepest valley of the Dyck word the *valley depth* or *v-depth* (for simplicity). For a Dyck number $d$ we denote the v-depth by $VD(d)$. In (8), v-depth is the smallest digit, i.e., $VD(d) = 0$.

The v-depth gives us the opportunity to refine the result of the DS-function. To get the successor of the Dyck number $d$ with binary expansion (7), we can fix the suffix 011111 (length 6, repunits 5) as Mersenne number $M_5$. For the new number $d'$, according to formula (6), we get a new suffix 100111, since the length of the zero fragment is $\lfloor 5/2 \rfloor = 2$ and repunits $5 - 2 = 3$. Well! But the result is not a Dyck number because some vertices of new path fall below the ground (bellow we replace four negative heights with '?'):

$$h(d') = 0?0???0121232321212321.$$

The situation can be corrected if we restore one unit in the suffix, i.e., 101111. In the end, in the original word, we simply swapped the trailing 0 and 1 on the right. This is how the DS-function algorithm often works. In binary expansion, the procedure for swapping two adjacent bits does not change the height of either the right bits or the left bits. Thus we get the following result:

$$DS(d) = 2{,}893{,}215 + 2^5 - 2^4 = 2{,}893{,}215 + 2^4 = 2{,}893{,}231$$
$$= 1011000100101101011111 \text{ (in base 2)};$$
$$h(DS(d)) = 2121012343454543434321.$$

In the last equality, only in one bit the height has changed (shown in large print); now it is a new valley with a height of 3. It is easy to see that replacing an arbitrary 0 with a 1 (or vice versa) changes the heights of subsequent bits by 2 up (or 2 down). So we come to the following proposition.

**Proposition 8.** For a given Dyck number, the DS-function is affected by two parameters: the v-depth and the repunit suffix length in the binary expansion.

For an arbitrary Dyck number, a rather simple plan for obtaining a successor arises. First, in binary expansion, we swap the right 0 and the adjacent 1 of the suffix. Then, if the v-depth is 2 or more, we bit by bit zero out the repunit suffix while the v-depth is greater than 1. Now we can proceed to the algorithm for calculating the DS-function on the non-zero Dyck numbers.

## 4  The algorithm

Let a non-zero Dyck number $d$ be given. It is necessary to calculate the successor of $d$, the Dyck number $d' = DS(d)$. Let's consider the DS-function algorithm step by step. The algorithm involves the following 8 steps:

*Step 1. Data initialization.* First you need to get the binary expansion of $d$. For further calculation we need three parameters: (a) the length of the binary code $n$, (b) the repunit suffix length $r$, and (c) the height of the deepest valley in the binary expansion, the v-depth $h = VD(d)$. A non-zero Dyck number is always odd, i.e., the binary expansion ends with 1, so $1 \leq r \leq n$. But we may not calculate the v-depth, since in some cases this parameter is not needed. In the first steps of the algorithm, we work as far as possible without the v-depth.

*Step 2*. Option $r = n$ (the binary expansion not contains zeros). The given Dyck number $d$ is the $n$-th Mersenne number. According to (4) $d' = d + 2^m$, $m = \lceil n/2 \rceil$. Output the value of $d'$ and finish the algorithm.

*Step 3*. Option $r = 1$ (there is only one unit in the repunit suffix). In the binary expansion we have the suffix 01 (the option $d = 1$, the 1-th Mersenne number, has already been processed in Step 2). As a result, we get a new Dyck number $d' = d + 2$ with the suffix 11 in the binary expansion. Let's pay attention to the fact that two Dyck numbers $d$ and $d'$ resemble a twin prime pair, moreover, they are often prime numbers, and there are a lot of such pairs (possibly an infinite number). Output $d'$ and finish the algorithm.

*Step 4*. Option $r = 2$ (in the repunit suffix 2 units). In this case we have the suffix 011 (the option $d = 3 = 11_2$, the 2-th Mersenne number, see Step 2). The result is a new binary suffix $101 = 11 + 10$ (in base 2). In fact, we just swapped the highest bit of the repunit suffix and the neighboring zero. We again get the equality $d' = d + 2$, and again we get an analogue of a twin prime pair with all the ensuing consequences. And the most interesting thing in this case is that the subsequent Dyck number is $d'' = d' + 2$ (see Step 3, suffix 01). As a result, we get a tuple of three odd numbers $(d, d + 2, d + 4)$. Let's call it *Dyck triplets*. Output $d'$ and finish the algorithm.

*Step 5*. Option $r = 3$ (there is three units in the repunit suffix). In the binary expansion we have the suffix 0111 (the option $d = 7 = 111_2$, the 3-th Mersenne number, has already been processed in Step 2). The result is a new binary suffix $1011 = 0111 + 1000 - 100$ (in base 2). As you can see, we have swapped the highest bit of the repunit suffix and the neighboring zero again. As a result, we get a new Dyck number $d' = d + 2^3 - 2^2 = d + 2^2 = d + 4$. Output the value of $d''$ and finish the algorithm.

*Step 6*. Option $3 < r < n$. Let's repeat the operation of the previous step, swap two bits in the binary expansion. We will swap the highest unit of the repunit suffix (the right peak in the corresponding Dyck path) and the neighboring zero on the left. This procedure reduces the height of the right peak by 2, and the peak becomes the right valley of the new path. The height of the other points of Dyck path does not change. As a result, we get a new Dyck number $d' = d + 2^r - 2^{r-1} = d + 2 \times 2^{r-1} - 2^{r-1} = d + 2^{r-1}$.

*Step 7*. Now we need a third parameter, v-depth $h$. In the new binary expansion of $d'$, we got a new right valley with height $r - 2$. And if it turns out that $h > r - 2$, then we need to adjust the v-depth: $h = VD(d') = r - 2$. Additionally, we check again if $h < 2$, then we issue $d'$ and finish the algorithm.

*Step 8*. In the binary expansion of $d'$, the length of the repunit suffix is $r' = r - 1 > 2$. We have a height reserve, $h \geq 2$; so we can decrease $d'$ if we start zeroing out the high-order bits of the repunit suffix. Recall that zeroing one digit of the suffix reduces $h$ by 2. In Step 6 in binary expansion, we swapped two adjacent bits. In particular, the $(r-1)$-th digit (the highest bit of the repunit suffix) has been zeroed. The result is the Dyck number $d' = d + 2^{r-1}$. Since we have a height reserve $h \geq 2$, let's also reset an additional $(r-2)$-th bit of the repunit suffix:

$$d" = d' - 2^{r-2} = d + 2^{r-1} - 2^{r-2} = d + 2 \times 2^{r-2} - 2^{r-2} = d + 2^{r-2}.$$

Obviously, if $h \geq 4$, then we will also set to zero the $(r-3)$-th bit of the suffix, and then we will get a new number $d''' = d + 2^{r-3}$ and so on. The number of zeroed bits depends on the value of $h$. Thus, we can organize a cycle in which the number of repetitions does not exceed $r/2$. At the end of the cycle, we output the resulting Dyck number and finish the algorithm.

Note that we can process the Mersenne number as an ordinary Dyck number if we insert one leading zero at the beginning of the binary expansion. The appendix shows the Python program for calculating a DS-function for an arbitrary Dyck number.

The above algorithm will be considered as a proof of the following theorem.

**Theorem 9.** Let $d$ be a Dyck number and let $r$ be the length of the repunit suffix in the binary expansion of $d$. Then

1. if $r = 0$, then $DS(d) = 1$;
2. if $r = 1$ or $r = 2$, then $DS(d) = d + 2$;
3. otherwise, $DS(d) = d + 2^{r-1-\lfloor h/2 \rfloor}$, where $h$ is the height of the deepest valley (v-depth) in the binary expansion after swapping the trailing 0 and 1 on the right.

Let's pay attention, if $h = r - 2$ (the valley before the repunit suffix is the deepest), we get

$$r - 1 - \lfloor (r-2)/2 \rfloor = r - 1 - (\lfloor r/2 \rfloor - 1) = r - \lfloor r/2 \rfloor = \lceil r/2 \rceil.$$

As a result, we obtain an analogue of the formula for Mersenne numbers, see (4):

(9) $$DS(d) = d + M_m + 1 = d + 2^m, \quad m = \lceil r/2 \rceil.$$

## 5 Conclusions and future work

This is the first paper devoted to Dyck numbers. We plan to publish several more papers in which we need to consider certain problems related to Dyck numbers. For example, it is interesting to know the length of the ranges in the OEIS A036991/A350346. And there is OEIS [A001405](A001405), in which the terms are quite suitable for the goal. Here is the beginning of A001405 (starting or offset from 0):

1, 1, 2, 3, 6, 10, 20, 35, 70, 126, 252, 462, 924, 1716, 3432, 6435, 12870 …

Let's show again the first ranges in A036991, but this time without the zero range. Still, we are working with real Dyck paths (Dyck words). Recall that many sources ignore empty paths, and now we need it.

$R_1 = \{1\}, |R_1| = 1 = A001405(0)$;
$R_2 = \{3\}, |R_2| = 1 = A001405(1)$;
$R_3 = \{5, 7\}, |R_3| = 2 = A001405(2)$;
$R_4 = \{11, 13, 15\}, |R_4| = 3 = A001405(3)$;
$R_5 = \{19, 21, 23, 27, 29, 31\}, |R_5| = 6 = A001405(4)$;

…………………..

$R_{16} = \{33023, 33151, \ldots, 65535\}$, $|R_{16}| = 6435 = $ A001405(15)   and so on.

It is easy to check the data on the A036991 b-file (in the b-file, the terms are numbered). In our case, the ordinal number of the range is determined by the length of the binary expansion of the terms, and we would be more comfortable if A001405 started at 1. We can formulate an appropriate conjecture.

**Conjecture 10.** The length of the *k*-th range in OEIS A036991 is A001405(*k*–1).

We cannot yet prove a similar theorem. The author will be grateful to interested readers who will take part in this project.

Acknowledgements. The author would like to thank Jörg Arndt  (Technische Hochschule Nürnberg),  Alois P. Heinz  (University of Applied Sciences Heilbronn), and Olena G. Kachko (Kharkiv National University of Radio Electronics) for discussion of the considered integer sequences.

# References

[1]  R. Stanley. *Catalan numbers*. Cambridge University Press, Cambridge, 2015.

[2]  Gennady Eremin. *Dynamics of balanced parentheses, lexicographic series and Dyck polynomials*, 2019. arXiv:1909.07675

[3]  Zoltan Kasa. *Generating and ranking of Dyck words*, Acta Universitatis Sapientiae, Informatica, **1**, 1 (2009), 109-118. arXiv:1002.2625

[4]  Neil J. A. Sloane, Ed. *The On-Line Encyclopedia of Integer Sequences*. https://oeis.org/, 2022.

[5]  Ralph L. Childress. *Recursive prime factorizations: Dyck words as Representations of numbers*, 2021. arXiv:2102.02777

[6]  Wolfram MathWorld, Mersenne Numbers, 2022.
https://mathworld.wolfram.com/MersenneNumber.html

[7]  Wolfram MathWorld, Repunit, 2022.
https://mathworld.wolfram.com/Repunit.html

[8]  R. Florez, T. Mansour, J. Ramirez, F. Velandia, D. Villamizar. *Restricted Dyck paths on valleys sequence*, 2021.  arXiv:2108.08299

Mentions the OEIS sequences A000225, A001405, A002275, A005408, A014486, A036991, A063171, A350346, A350577.

Gzhel State University, Moscow, 140155, Russia
http://www.art-gzhel.ru/

**Appendix.** The Python program calculates a DS-function for any Dyck number.

```python
def repunit(n):     # calculate the length of the repunit suffix
    rep = 0
    for bit in bin(n)[:1:-1]:
        if bit == '1': rep += 1;
        else break
    return rep

def val_deep(n):    # the deepest valley
    height, key, val = 0, 1, 999
    for bit in bin(n)[:1:-1]:
        if bit == '1':
            if key == 0 and height < val: val = height
            height += 1; key = 1
        else: height -= 1; key = 0
    return val

def dyck_succ(n):   # dyck successor
    if n == 0: return 1
    ru = repunit(n)       # the length of repunit suffix
    if ru < 3: return n + 2  # triplets
    val = val_deep(n)     # the deepest valley in path
    if val > ru-2: val = ru - 2
    return n + 2**(ru - 1 - val//2)

d = 65535  # slightly expand b-file in OEIS A036991
for n in range(13496,14001):
    print(n, d, bin(n))  # binary expansion for analysis
    d = dyck_succ(d)
```